\documentclass[reqno]{amsart}
\usepackage{latexsym}
\usepackage[latin1]{inputenc}
\usepackage{amsfonts}
\usepackage{enumerate}

\newtheorem{lema}{Lemma}
\newtheorem{teo}{Theorem}

\newcommand{\iz}[1]{{\flushleft \bf #1}}
\newcommand{\cal}{\mathcal}
\title[Existence and Stability of Almost Periodic Solutions of DEPCAG]{Existence and Stability of Almost Periodic Solutions of Differential Equations with Generalized Piecewise Constant Argument}
\author[S. Castillo \& M. Pinto]{Samuel Castillo {\bf \textdagger}
and Manuel Pinto {\bf \textdaggerdbl}}
\address{Samuel Castillo. Departamento de Matem\'atica.
Facultad de Ciencias. Universidad del B\'{\i}o-B\'{\i}o. Casilla
5-C. Concepci\'on. Chile.}\email{scastill@ubiobio.cl}
\thanks{{\bf \textdagger} Supported by DIUBB 110908 2/R,  FONDECYT 1080034}\thanks{{\bf \textdaggerdbl} Supported by
FONDECYT 1080034, FONDECYT  1120709 and DGI MATH UNAP 2009}
\address{Manuel Pinto.
Departamento de Matem\'atica. Facultad de Ciencias. Universidad de
Chile. Casilla 653. Santiago. Chile.}\email{pintoj@uchile.cl}
\date{}
\begin{document}
\begin{abstract}
This work deals with the existence of an almost periodic solution
for certain kind of differential equations with generalized
piecewise constant argument, almost periodic coefficients which are
seen as a perturbation of a linear equation of that kind satisfying
an exponential dichotomy on a difference equation. The stability of
that solution in a semi-axis studied.
\end{abstract}
\maketitle \iz{2010 AMS Subject Class:} 34A38, 34C27, 34D09, 34D20.
\iz{Key words:} Differential Equations, Piecewise constant argument,
Almost Periodic Solutions, Exponential Dichotomy, Stability of
Solutions.

\section{Introduction}
Let $\mathbb{N}$, $\mathbb{Z}$, $\mathbb{R}$, $\mathbb{C}$ be the
sets of natural, integer, real and complex numbers, respectively.
Denote by $|\cdot|$ the Euclidean norm for every finite dimensional
space on $\mathbb{R}$.

Fix a real valued sequence $(t_n)_{n=-\infty}^{+\infty}$,  such that
$t_n<t_{n+1}$ and $t_n \to \pm\infty$ as $n \to \pm\infty$. For $p
\in \mathbb{Z}$, let $\gamma^{p}: \mathbb{R} \to \mathbb{R}$ be
functions such that $\gamma^{p}/J_n=t_{n-p}$ for all $n \in
\mathbb{Z}$, where $J_n=[t_n,t_{n+1}[$, for all $n \in \mathbb{Z}$.

We are interested in the existence of almost periodic solution of
the following linear {\it Differential Equations with Piecewise
Constant Argument Generalized} (DEPCAG)
\begin{equation}
\label{depca1} y'(t)=A(t)y(t)+B(t)y(\gamma^{0}(t))+f(t),\;t\in
\mathbb{R}
\end{equation}
and
\begin{equation}
\label{depca2}
y'(t)=A(t)y(t)+B(t)y(\gamma^{0}(t))+F(t,y_{\gamma}(t)),\;t\in
\mathbb{R},
\end{equation}
where
\begin{equation}
\label{ygama}
y_{\gamma}(t)=(y(\gamma^{p_1}(t)),y(\gamma^{p_2}(t)),...,y(\gamma^{p_{\ell}}(t))),
\end{equation}
where $p_1,p_2,\ldots,p_\ell \in \mathbb{N}\cup \{0\}$.
DEPCAGs (\ref{depca1}) and (\ref{depca2}) are seen as perturbation of the linear DEPCAG
\begin{equation}
\label{homo2} z'(t)=A(t)z(t)+B(t)z(\gamma^0(t)),
\end{equation}
where $A,B: \mathbb{R} \to {\cal M}_{q}(\mathbb{C})$ and $f:\mathbb{R} \to \mathbb{C}^q$ are locally integrable functions and $F:\mathbb{R} \times W \subseteq \mathbb{R} \times(\mathbb{C}^q)^{\ell}\to \mathbb{C}^q$ is a continuous function.

For our study, the following additional assumptions are made.

\begin{enumerate}
\item[(H1)] $A$ and $B$  are almost periodic functions.
\item[(H2)] $\displaystyle\left(t_{n}^{(k)}\right)_{n=-\infty}^{+\infty}$, where $
t_{n}^{(k)}=t_{n+k}-t_n$ for all $k \in \mathbb{Z}$, is
equipotentially almost periodic for all $k \in \mathbb{Z}$.
\item[(H3)] (H2) holds and for all $\varepsilon>0$,
\[
T(f,\varepsilon)=\left\{\tau \in \mathbb{R}:|f(t+\tau)-f(t)| \leq
\varepsilon,\;\forall
t\in\mathbb{R}-\left(\bigcup_{n\in\mathbb{Z}}]t_n-\varepsilon,t_n+\varepsilon[\right)\right\}
\]
is relatively dense and  there is $\delta_{\varepsilon}>0$ such that
$|f(t'+\tau')-f(t')| \leq \varepsilon$ if  $\tau' \in
\mathbb{R}:|\tau'|\leq \delta_{\varepsilon}$ and $t',t'+\tau'$ is in
some of the intervals $[t_n,t_{n+1}]$.
\item[(H4)] $F$ is uniformly almost periodic on $W$ and there is $L>0$  such that
\begin{equation}
\label{lp} |F(t,x_1,...,x_{\ell})-F(t,y_1,...,y_{\ell})| \leq
L\sum_{j=1}^{\ell} |x_j-y_j|,
\end{equation}
for all $t \in \mathbb{R}$ and
$(x_1,...,x_{\ell}),(y_1,...,y_{\ell}) \in W$.
\end{enumerate}
A kind of exponential dichotomy is imposed on a part of the linear
DEPCAG \ref{homo2}, which will be made explicit in the following
section.

This work is motivated by the results in Fink \cite[Theorem 7.7,
Theorem 8.1 and Theorem 11.31]{F5}. Some extensions for piecewise
constant argument can be found in \cite{A02,YH16,Y}. Existence of
almost periodic solutions for the impulsive case can be found in
\cite{PR,SP} Our focus is to see the almost periodic solutions for
DEPCAGs (\ref{depca1}) and (\ref{depca2}) in terms of the solutions
of the difference equation from the Cauchy Operator of the linear
part (\ref{homo2}), on the points $t_n$ for all $n \in \mathbb{N}$,
in the style of \cite{YH16}. This work is different to Akhmet work
\cite{A02} where an exponential dichotomy on an ordinary
differential system is considered. This work is different to the
works on Hong-Yuan \cite{YH16} and Yuan \cite{Y} since a more
general $y_{\gamma}$ is considered.

Let $X$ be a fundamental matrix of the linear homogeneous system
\begin{equation}
\label{homo} x'=A(t)x
\end{equation}
and $X(t,s)=X(t)X(s)^{-1}$. Now we follows \cite{A03} to say what is the Cauchy matrix for (\ref{homo2}).

For $n \in \mathbb{Z}$ and $t\in J_n$
such that $t \geq s$, let  $Z_n(t)=X(t,t_n){\cal J}_n(t)$, where $\displaystyle {\cal
J}_n(t)=I+\int_{t_n}^t X(t_n,u)B(u)du$ and assume that
\begin{equation}
\label{invj} {\cal J}_{n}(t)\;\mbox{is invertible, for all}\;n \in \mathbb{Z}\;\mbox{and}\;t \in [t_n,t_{n+1}].
\end{equation}

Let
\begin{equation}
\label{oc} H(n)=Z_n(t_{n+1}),
\end{equation}
for all $n \in \mathbb{Z}$. For $\tau \in \mathbb{R}$, let $k(\tau)
\in \mathbb{Z}$ such that $\tau \in J_{k(\tau)}$. Consider $t>s$
such that $k(t)>k(s)$. Then, it is defined
\begin{equation}
\label{productoria} Z(t,s)=Z_{k(t)}(t)\left[
H(k(t)-1)H(k(t)-2)\cdots
H(k(s)+1)\right]H(k(s))^{-1}Z_{k(s)}(s)^{-1}.
\end{equation}
If $t \leq s$, by condition (\ref{invj}), $Z(t,s)=Z(s,t)^{-1}$ is
well defined. So, $Z(t,s)$ is the Cauchy matrix for (\ref{homo2}).
(see \cite{A,A02,Pi2010JDEA,ShWi83,Wi83,Wi93}).

Consider the difference equation
\begin{equation}
\label{disdich} \phi({n+1})=H(n)\phi(n).
\end{equation}
Notice that if $z:\mathbb{R} \to \mathbb{C}$, then $\phi(n)=z(t_n)$
is a solution of (\ref{disdich}) if $z$ is a solution of
(\ref{homo2}).

It will  be proved that $H=(H(n))_{n=-\infty}^{+\infty}$ in
(\ref{oc}) is almost periodic and that the sequence
$h=(h(n))_{n=-\infty}^{+\infty}$, defined by
\begin{equation}
\label{oh} h(n)=\int_{t_n}^{t_{n+1}} X(t_{n+1},u)f(u)du,
\end{equation}
for all $n \in \mathbb{Z}$, is almost periodic. Based on the
exponential dichotomy of (\ref{disdich}) and the almost periodicity
of $H$ and $h$, it will be proved that the bounded solution $c$ of
the discrete system
\begin{equation}
\label{disc} c(n+1)=H(n)c(n)+h(n),
\end{equation}
is almost periodic and the correspondence $h \mapsto c$ is Lipschitz continuous. Then it will be proved that the inhomogeneous
linear DEPCAG (\ref{depca1}) has an analogous almost periodic solution. The dependence of the almost periodic solution can be seen in terms
of the almost periodic solution of the discrete part for
(\ref{depca1}) and (\ref{depca2}), the linear continuous dependence
of the almost periodic solution $y$ of (\ref{depca1}) in terms of
$f$ and the same kind of dependence of $c$ of the almost periodic
solution of (\ref{disc}) in terms of $h$.

 By
assuming that $L$ in (\ref{lp}) is small enough, an almost periodic
solution for the DEPCAG (\ref{depca2}) is obtained in terms of the
solution of a difference equation. Finally, it will be proved that
the almost periodic solution of the DEPCAG (\ref{depca2}) is
exponentially stable as $t \to +\infty$ with respect the solutions
of (\ref{depca2}) for $t \geq 0$.
The exponential stability is proved by using a
Gronwall inequality on the mentioned difference equation.

The present work is organized as follows: Section 2 provides the
main definitions, assumptions and facts that will be used. In the
Section 3, the existence of almost periodic solutions for the DEPCAG
(\ref{depca1}) is studied. In Section 4, that study is extended for
the DEPCAG (\ref{depca2}) and deals with asymptotic stability
for the DEPCAG (\ref{depca2}) as $t \to +\infty$. An
example is given in the last section.

\section{Preliminaries}

\noindent (H6) Assume that (\ref{disdich}) has an {\it exponential
dichotomy.}

The recent assumption is equivalent to assume that there is a
projection $\Pi:\mathbb{C}^q \to \mathbb{C}^q$ and positive
constants $\rho, K$ with $\rho<1$ such that
\begin{equation}
\label{edic} |{\cal G}(n,k)| \leq K\rho^{\pm(n-k)},
\end{equation}
for all $n,k \in \mathbb{Z}:\pm(n-k) \leq 0$, where
\begin{equation}
\label{green} {\cal G}(n,k)=\left\{\begin{array}{rcl}
\Phi(n)\Pi\Phi(k+1)^{-1},&\mbox{if}& n>k\\
\\
-\Phi(n)(I-\Pi)\Phi(k+1)^{-1},&\mbox{if}& n \leq k
\end{array}\right.
\end{equation}
and $\Phi$ is a fundamental matrix for the system (\ref{disdich}).
In particular it will be said that system (\ref{disdich}) {\it
exponentially stable} as $n \to +\infty$ if it has an exponential
dichotomy with $\Pi=I$.

The recent dichotomy definition has been adapted from that given by
Papashinopulos \cite{Papas94} for (\ref{homo2}) when
$\gamma=[\cdot]$. It is an exponential dichotomy for (\ref{disdich})
which is not obvious to be extended for (\ref{homo2}) in terms of
$Z(t,s)$ except for cases where the projection for exponential
dichotomy commutes with $A(t)$ and $B(t)$. Authors has not
found
any reference containing a definition of exponential dichotomy for
(\ref{homo2}).

We start with some classical notions.

It is said that $x$ is a solution of a DEPCAG
\begin{equation}
\label{0}
x'(t)=\tilde{f}(t,x_{\gamma}(t)),
\end{equation}
where $x_{\gamma}$ is defined in (\ref{ygama}), if
\begin{enumerate}[(a)]
\item $x$ is continuous on $\mathbb{R}$;
\item the derivative $x'$ of $x$ exists except possibly at the points $t=t_n$ with $n \in \mathbb{Z}$, where every one-sided derivative exist;
\item $x$ is a solution of the DEPCAG (\ref{0}) except possibly at the points $t=t_n$ with $n \in \mathbb{Z}$.
\end{enumerate}
 If $\mathbb{E}$ is a
finite dimensional space on $\mathbb{R}$, $D \subseteq \mathbb{R}$
and $g:D \to \mathbb{E}$, then $\displaystyle |g|_{\infty}=\sup_{t
\in D} |g(t)|$.  A set $E \subseteq \mathbb{R}$ is called {\it
relatively dense} if there exists a positive real number $l$ such
that $E \cap [m,m+l] \neq \phi$ for all $m \in \mathbb{R}$. For
$\mathbb{A}\subseteq \mathbb{R}$ an additive group and
$\left(\mathbb{E},|\cdot|\right)$ a finite dimensional linear space
$g:\mathbb{A} \to \mathbb{E}$ is called {\it almost periodic} if it
is continuous the set of translations $T(g,\varepsilon)$, defined by
the set of all $\tau \in \mathbb{A}$ such that $|g(t+\tau)-g(t)|
\leq \varepsilon$ for all $t\in \mathbb{A}$, is relatively dense for
all $\varepsilon>0$ (see \cite[Definition 1.10]{F5}). There will be
considered the cases $\mathbb{A}=\mathbb{R}$ (almost periodic
functions) and $\mathbb{A}=\mathbb{Z}$ (almost periodic sequences).
We can notice by following \cite[page 201]{SP} that (H3),  is a
definition of almost periodicity for {\it piecewise} continuous
functions. An alternative definition of almost periodicity {\it for
continuous functions} was given by Salomon Bochner \cite{Bochner}
(see Fink \cite[page 14]{F5} for more detailed reference): A
function $f$ is almost periodic if every sequence $\left(f(t_n +
t)\right)_{n=1}^{+\infty}$ of translations of $f$ has a subsequence
that converges uniformly for $t \in \mathbb{R}$. A function
$F:\mathbb{R} \times W \subseteq \mathbb{R} \times \mathbb{E}\to
\mathbb{E}^q$ is uniformly almost periodic on W, if the set
$T(F,\varepsilon,W)$ which denotes the set of all $\tau \in
\mathbb{R}$ such that $|F(t+\tau,w)-F(t,w)| \leq \varepsilon$ for
all $(t,w) \in\mathbb{R} \times W$, is relatively dense for every
$\varepsilon>0$.

Next, some notations are given.
 Let ${\cal A}{\cal
P}(\mathbb{A},\mathbb{E})$ be the set of the almost periodic
functions from $\mathbb{A}$ into $\mathbb{E}$. The set  $\left({\cal
A}{\cal P}(\mathbb{A},\mathbb{C}^q),|\cdot|_{\infty}\right)$ is a
Banach space.

We say that $\left(t_n^{(k)}\right)_{n=-\infty}^{+\infty}$ is
{\it equipotentially almost periodic}, for all $k \in \mathbb{Z}$ if the set
\[
 \bigcap_{k \in \mathbb{N}} \left\{T \in \mathbb{Z}:\left|t_{T+n}^{(k)}-t_{n}^{(k)}\right| \leq \varepsilon,\;\mbox{for all}\; n \in \mathbb{Z}\right\}
\]
is relatively dense for all $\varepsilon>0$.

Since $A$, $B$  are almost periodic, $A$, $B$  are
bounded. Since $\left(t_n^{(k)}\right)_{n=-\infty}^{+\infty}$
is equipotentially almost periodic for all $k \in  \mathbb{Z}$, every sequence
$\left(t_n^{(k)}\right)_{n=-\infty}^{+\infty}$ is almost periodic
for all $k \in \mathbb{Z}$. So, the sequences
$\left(t_n^{(k)}\right)_{n=-\infty}^{+\infty}$ are bounded for all
$k \in \mathbb{Z}$ (see \cite[Theorem 67]{SP}) and there exists the
positive real number
\begin{equation}
\label{theta} \theta=\sup_{n \in \mathbb{Z}} (t_{n+1}-t_n).
\end{equation}

Since
\[
|Z(t,s)| \leq
e^{|A|_{\infty}(t_{n+1}-t_n)}\left(1+e^{|A|_{\infty}(t_{n+1}-t_n)}|B|_{\infty}(t_{n+1}-t_n)\right),
\]
for all $t,s \in J_n$. So, $Z(t,s)$ is bounded. By following \cite{A03,Pi2010JDEA}, we have that $y:\mathbb{R}
\to \mathbb{C}^q$ given by
\begin{equation}
\label{pre-22}
\begin{array}{rcl}
y(t)&=&Z_{k(t)}(t)\times
\left(\sum_{k=-\infty}^{+\infty}{\cal G}(k(t),k)\int_{t_k}^{t_{k+1}} X(t_{k+1},u)f(u)du\right)\\
\\
&+&\int_{\gamma^0(t)}^t X(t,u)f(u)du,
\end{array}
\end{equation}
where $t \in \mathbb{R}$, will be the unique bounded solution of
(\ref{depca1}) which satisfies (H3) (see Theorem 2 below). Moreover,
by taking limits $t \to \gamma^0(t)^+$ and $t \to \gamma^0(t)^-$, we
obtain that $y$ is continuous  en every $t_n$ and therefore $y$ is
almost periodic.

For $\varepsilon>0$, let $\Gamma_{\varepsilon}$ be the set of $r \in
\mathbb{R}$ such that there is $k \in \mathbb{Z}$ with
\begin{equation}
\label{disc-t} \sup_{n \in \mathbb{Z}}\left|t_{n}^{(k)}-r\right|
\leq \varepsilon.
\end{equation}
Denote by $P_{r}(\varepsilon)$ the set of all $k \in \mathbb{Z}$
satisfying (\ref{disc-t}). Let
\[
P_{\varepsilon}=\bigcup_{r \in
\Gamma_{\varepsilon}}P_{r}(\varepsilon).
\]

We will need the following lemmas.
\begin{lema} \cite[Lemma 23]{SP}
\label{leap1} Assume that (H2) holds. Let $\varepsilon>0$, $\Gamma \subseteq
\Gamma_{\varepsilon}$, $\Gamma \neq \phi$ and $\displaystyle P
\subseteq \bigcup_{r \in \Gamma} P_{r}(\varepsilon)$ be such that $P
\cap P_{r}(\varepsilon) \neq \phi$ for all $r \in \Gamma$. Then the
set $\Gamma$ is relatively dense if and only if $P$ is relatively
dense.
\end{lema}

\begin{lema} \cite[Lemma 25]{SP}
\label{leap2} The following statements are equivalent.
\begin{enumerate}[(a)]
\item (H2) holds;
\item The set $P_{\varepsilon}$ is relatively dense for any $\varepsilon>0$;
\item The set $\Gamma_{\varepsilon}$ is relatively dense for any $\varepsilon>0$.
\end{enumerate}
\end{lema}

\begin{lema} \cite[Lemma 29]{SP}
\label{rdg} Assume that $f$ satisfies (H3). Then
$\Gamma_{\varepsilon} \cap T(f,\varepsilon)$ is relatively dense.
\end{lema}

By mean standard arguments, it can be proved the following result.
\begin{lema}
\label{master}
\begin{enumerate}[(a)]
\item If $f_1, f_2$ are functions satisfying (H3), then given $\varepsilon>0$,\\ $\Gamma_{\varepsilon} \cap T(f_1,\varepsilon)\cap
T(f_2,\varepsilon)$ is relatively dense.
\item If $(g_1(n))_{n=-\infty}^{+\infty}$ and
$(g_2(n))_{n=-\infty}^{+\infty}$ are almost periodic solutions, then
given $\varepsilon>0$, $P_{\varepsilon} \cap T(g_1,\varepsilon)\cap
T(g_2,\varepsilon)$ is relatively dense.
\end{enumerate}
\end{lema}

For the following results, remind that $q$ is the dimension of
DEPCAG (\ref{homo2}). Notice that they depends only on the
assumptions (H1) and (H3).

\begin{lema}
\label{conti} Consider $\theta$ defined in (\ref{theta}). Let
$K_0=\exp(|A|_{\infty}\theta)$, $\displaystyle K_1= \sup_{n \in
\mathbb{Z}}e^{\left(|A|_{\infty}\left|t_{n+1}^{(p)}-\tau\right|\right)}$
and $\displaystyle K_2=K_0K_1$. Then
\begin{enumerate}[(a)]
\item $|X(t,s)| \leq \sqrt{q}K_0$, for all $t,s \in \mathbb{R}$ such that $|s-t|\leq \theta$;
\item If $\tau>0$, $p \in \mathbb{N}$ and $u \in [t_n,t_{n+1}]$  then
\[
\begin{array}{rcl}
|X(t_{n+p+1},u+\tau)-X(t_{n+1},u)| &\leq& \sqrt{q}\cdot\left[K_1|A|_{\infty}\left|t_n^{(p)}-\tau\right|\right.\\
\\
&+&\left.K_2|A(\cdot+\tau)-A(\cdot)|_{\infty}|t_{n+1}-t_n|\right]\\
\\
&\times&\exp\left(|A|_{\infty}(t_{n+1}-t_n)\right);\\
\\
\end{array}
\]

\item If $\tau>0$, $p \in \mathbb{N}$ and $t \in [t_n,t_{n+1}]$  then
\[
\begin{array}{rcl}
|X(t+\tau,t_{n+p})-X(t,t_{n})| &\leq& \sqrt{q}\cdot\left[K_1\left|t_n^{(p)}-\tau\right|\right.\\
\\
&+&\left.K_2|A(\cdot+\tau)-A(\cdot)|_{\infty}\right]\left|t_{n+1}^{(p)}-\tau-(t_{n+1}-t_n)\right|\\
\\
&\times&\exp\left(|A|_{\infty}\left(\left|t_n^{(p)}-\tau\right|+\theta\right)\right);\\
\\
\end{array}
\]

\item If $\tau>0$ and $t,s \in \mathbb{R}: |t-s| \leq \theta$ then
\[
|X(t+\tau,s+\tau)-X(t,s)| \leq \sqrt{q}K_0
|A(\cdot+\tau)-A(\cdot)|_{\infty};
\]
\item If $\tau>0$, $p \in \mathbb{N}$ and $u \in [t_n,t_{n+1}]$ then
\[
\begin{array}{rcl}
|X(t_{n+p+1},t_{n+p})-X(t_{n+1},t_{n})|&\leq&K_2|X(u+\tau,t_{n+p})-X(u,t_{n})|\\
\\
&+&\sqrt{q}K_0|X(t_{n+p+1},u+\tau)-X(t_{n+p},u)|.
\end{array}
\]
\end{enumerate}
\end{lema}
\iz{Proof:} Part (a) follows immediately. To prove (b), assume
without loss of generality that $t_{n+p+1}-\tau \geq t_{n+1}$.
Notice that for $u \in [t_n,t_{n+1}]$,
\[
\begin{array}{rcl}
\Delta_n(u)&=&\int_{t_{n+1}}^{t_{n+p+1}-\tau} X(t_{n+p+1},\xi+\tau)A(\xi+\tau)d\xi\\
\\
&+&\int_{u}^{t_{n+1}} X(t_{n+p+1},\xi+\tau)[A(\xi+\tau)-A(\xi)]d\xi\\
\\
&+&\int_{u}^{t_{n+1}} \Delta_n(\xi)A(\xi)d\xi,\\
\\
\end{array}
\]
where $\Delta_n(u)=X(t_{n+p+1},u+\tau)-X(t_{n+1},u)$. Then
\[
\begin{array}{rcl}
|\Delta_n(u)|&\leq&\int_{t_{n+1}}^{t_{n+p+1}-\tau} |X(t_{n+p+1},\xi+\tau)||A(\xi+\tau)|d\xi\\
\\
&+&\int_{t_n}^{t_{n+1}} |X(t_{n+p+1},\xi+\tau)||A(\xi+\tau)-A(\xi)|d\xi\\
\\
&+&\int_{u}^{t_{n+1}} |\Delta_n(\xi)||A(\xi)|d\xi.\\
\\
\end{array}
\]

So, by Gronwall's inequality the result is obtained.

Similarly, assume without loss of generality that $t_{n+1} \geq
t_{n+p}-\tau$. If $\Delta_n^*(t)=X(t+\tau,t_{n+p})-X(t,t_n)$, then
\[
\begin{array}{rcl}
|\Delta_n^*(t)|& \leq &\left|\int_{t_{n+p}-\tau}^{t_n} |A(\xi)||X(\xi+\tau,t_{n+p})|d\xi\right|\\
\\
&+&\int_{t_{n+p}-\tau}^{t_{n+1}} |A(\xi+\tau)-A(\xi)||X(\xi+\tau,t_{n+p})|d\xi\\
\\
&+&\int_{t_n}^{t} |A(\xi)||\Delta_n^*(\xi)|d\xi,\\
\\
\end{array}
\]
for $t \in [t_{n},t_{n+1}]$. So, by Gronwall's inequality, (c) is
obtained. To prove part (d), proceed as in the proof of
\cite[Proposition 8]{YH16}. To prove (e), notice that
\[
\begin{array}{rcl}
X(t_{n+p+1},t_{n+p})-X(t_{n+1},t_{n})&=&X(t_{n+p+1},u+\tau)[X(u+\tau,t_{n+p})-X(u,t_{n})]\\
\\
&+&[X(t_{n+p+1},u+\tau)-X(t_{n+p},u)]X(u,t_n)
\end{array}
\]
and apply the previous results. $\Box$


By Lemma \ref{conti}, the following result is obtained.

\begin{lema}
\label{zoom} Consider $\theta$ defined in (\ref{theta}). Let
$\varepsilon>0$, $\tau \in \Gamma_{\varepsilon}\cap
T(A,\varepsilon)$ and $p \in P_{\tau}(\varepsilon)$. Then there is
$K' >0$ such that for all $n \in \mathbb{Z}$
\begin{enumerate}[(a)]
\item $|X(t_{n+p+1},u+\tau)-X(t_{n+1},u)| \leq K'\varepsilon$,
for all $u \in [t_n,t_{n+1}]$;
\item $|X(t+\tau,t_{n+p})-X(t,t_n)| \leq K'\varepsilon$,
for all $t \in [t_n,t_{n+1}]$;
\item $|X(t+\tau,s+\tau)-X(t,s)| \leq K'\varepsilon$, for all $s,t \in \mathbb{R}:|t-s|\leq \theta$;
\item $|X(t_{n+p+1},t_{n+p})-X(t_{n+1},t_{n})|\leq K'\varepsilon$.
\end{enumerate}
\end{lema}

Notice that the existence of $p \in P_{\tau}(\varepsilon)$ is given
by Lemma \ref{leap2} and the existence of $\tau \in
\Gamma_{\varepsilon}\cap T(A,\varepsilon)$ is given by Lemma
\ref{rdg}.

\begin{lema}
\label{apH} The sequence $H=(H(n))_{n=-\infty}^{+\infty}$ given by
(\ref{oc}) and the sequence $h=(h(n))_{n=-\infty}^{+\infty}$ given
by (\ref{oh}) are almost periodic.
\end{lema}
\iz{Proof:} Firstly, notice that $\displaystyle
H(n)=X(t_{n+1},t_n)+\psi(n)$, for all $n \in \mathbb{Z}$, where
\[
\psi(n)=\int_{t_n}^{t_{n+1}} X(t_{n+1},u)B(u)du.
\]

From Lemma \ref{zoom} (d), it is not hard to see that
$\left(X(t_{n+1},t_n)\right)_{n=-\infty}^{+\infty}$ is almost
periodic. $\psi$ is also almost periodic. In fact, let
$\varepsilon>0$.  From Lemma \ref{master},  $
\Gamma=T(A,\varepsilon)\cap T(B, \varepsilon)\cap
\Gamma_{\varepsilon}$ is relatively dense. Let $\displaystyle p \in
P=\bigcup_{\tau \in \Gamma} P_{\tau}(\varepsilon)$, so there is
$\tau \in \Gamma$ such that $p \in P_{\tau}(\varepsilon)$. Then, for
all $n \in \mathbb{Z}$ it is obtained
\[
\begin{array}{rcl}
\psi(n+p)-\psi(n)& =&
\int_{t_{n+p}}^{t_{n+p+1}}X(t_{n+p+1},u)B(u)du-
\int_{t_{n}}^{t_{n+1}}X(t_{n+1},u)B(u)du\\
\\
&=&\int_{t_{n+p}}^{t_{n+p+1}}X(t_{n+p+1},u)B(u)du
-\int_{t_{n}+\tau}^{t_{n+p+1}}X(t_{n+p+1},u)B(u)du\\
\\
&+& \int_{t_{n}+\tau}^{t_{n+p+1}}X(t_{n+p+1},u)B(u)du
-\int_{t_{n}}^{t_{n+1}}X(t_{n+p+1},u+\tau)B(u+\tau)du\\
\\
&+& \int_{t_{n}}^{t_{n+1}}X(t_{n+p+1},u+\tau)B(u+\tau)du
-\int_{t_{n}}^{t_{n+1}}X(t_{n+1},u)B(u)du,\\
\\
& =&
\int_{t_{n+p}}^{t_{n}+\tau}X(t_{n+p+1},u)B(u)du+
\int_{t_{n+1}+\tau}^{t_{n+p+1}}X(t_{n+p+1},u)B(u)du
\\
\\
&+&
\int_{t_{n}}^{t_{n+1}}\left[X(t_{n+p+1},u+\tau)B(u+\tau)-X(t_{n+1},u)B(u)\right]du.\\
\\
\end{array}
\]
By Lemmas \ref{conti} and \ref{zoom}, there are positive constants
$M$ and $K'$ such that
\[
\begin{array}{rcl}
|\psi(n+p)-\psi(n)|& \leq &
\left|t_n^{(p)}-\tau\right|M+
\left|t_{n+1}^{(p)}-\tau\right|M
+
K'\varepsilon,\\
\\
&\leq&[2M+K']\varepsilon\\
\\
\end{array}
\]
for all $n \in \mathbb{Z}$. So, $p \in T(\psi,[2M+K']\varepsilon)$.
Since $p$ was taken arbitrarily in $P$, $P \subseteq
T(\psi,[2M+K']\varepsilon)$. By Lemma \ref{leap1}, $P$ is relatively
dense. So, $T(\psi,[2M+K']\varepsilon)$ is relatively dense. Since
$\varepsilon>0$ is arbitrary, $\psi$ is almost periodic. Therefore,
$H=(H(n))_{n=-\infty}^{+\infty}$ is almost periodic.

In the similar way, $h$ is almost periodic. $\Box$

\section{Inhomogeneous Linear DEPCAG}
To study the existence of an almost periodic solution of  the DEPCAG
(\ref{depca1}), it is reminded that $\displaystyle f \in {\cal
A}{\cal P}(\mathbb{R},\mathbb{C}^q)$.

By mean of the constant variation formula \cite{A03,Pi2010JDEA},
\begin{equation}
\label{vpn} y(t)=Z(t,k(t)) c(k(t))+\int_{\gamma^0(t)}^t
X(t,u)f(u)du,
\end{equation}
is obtained, for all $t\in\mathbb{R}$, where $c$ is solution of the discrete
system (\ref{disc}). By taking $t \to t_{n+1}^-$, it is obtained a
solution $y$ for (\ref{depca1}) such that $y(t_n)=c(n)$ for all $n
\in \mathbb{Z}$. It will be proved that $y$ is almost periodic.

If $c$ is the bounded solution of equation (\ref{disc}) then
\begin{equation}
\label{fpo} c(n)=\sum_{k=-\infty}^{+\infty}{\cal G}(n,k)h(k),
\end{equation}
where the Green matrix ${\cal G}(n,k)$ is given by (\ref{green}) and
$h$ is given by (\ref{oh}).

From (\ref{vpn}) and (\ref{fpo}), $y$ is the bounded solution of
(\ref{depca1}) and satisfies (\ref{pre-22}).
 This relation shows $y$ as
a bounded linear function of $f$.

 By using the equivalent definition of almost periodicity due to S. Bochner, two important facts are obtained.
\begin{lema}(\cite[Proposition 7]{YH16} and \cite{Z})
\label{p7} A sequence $x=(x(n))_{n=-\infty}^{+\infty}$ is almost
periodic if and only if for any integer sequences
$(k_j')_{j=1}^{+\infty}$ and $(\ell_j')_{j=1}^{+\infty}$ there are
subsequences $k=(k_j)_{j=1}^{+\infty}$ and
$\ell=(\ell_j)_{j=1}^{+\infty}$ of $(k_j')_{n=1}^{+\infty}$ and
$(\ell_j')_{n=1}^{+\infty}$ respectively, such that
\[
T_kT_{\ell}x=T_{k+\ell}x,
\]
uniformly on $\mathbb{Z}$, where
$k+\ell=(k_j+\ell_j)_{j=1}^{+\infty}$, $\displaystyle
T_mx(n)=\lim_{j \to +\infty} x(n+m_j)$ and $m=(m_j)_{j=1}^{+\infty}
\in \{k,\ell,k+\ell\}$, for all $n \in \mathbb{Z}$.
\end{lema}

\begin{teo}
\label{hch}  Assume that hypotheses (H1), (H3) and (H6) are
satisfied. If  $c$ is given by (\ref{fpo}), then $c$ is the unique
almost periodic solution of the linear inhomogeneous difference
system (\ref{disc}). Moreover,
\begin{equation}
\label{estim-c} |c|_{\infty} \leq \frac{2K}{1-\rho}|h|_{\infty}.
\end{equation}
\end{teo}


\iz{Proof:} By Lemmas \ref{apH}, \ref{hap} and \ref{p7}, for any
integer sequences $(k_j')_{j=1}^{+\infty}$ and
$(\ell_j')_{j=1}^{+\infty}$ there are subsequences
$k=(k_j)_{j=1}^{+\infty}$ and $\ell=(\ell_j)_{j=1}^{+\infty}$ of
$(k_j')_{n=1}^{+\infty}$ and $(\ell_j')_{n=1}^{+\infty}$
respectively, such that $T_{k+\ell}H=T_kT_{\ell}H$ and
$T_{k+\ell}h=T_kT_{\ell}h$, uniformly on $\mathbb{Z}$.

Now, notice that $c$ given by (\ref{fpo}) is the only solution of
(\ref{disc}) which is bounded. Moreover, $z=T_{k+\ell}c$ and
$z=T_kT_{\ell}c$ are bounded solutions of
\[
\begin{array}{rcl}
z(n+1)&=&T_{k+\ell}H(n)z(n)+T_{k+\ell}h(n),\\
\\
z(n+1)&=&T_kT_{\ell}H(n)z(n)+T_kT_{\ell}h(n),\\
\\
\end{array}
\]
respectively. By uniqueness $T_{k+\ell}c=T_kT_{\ell}c$. So,
$c=(c(n))_{n=-\infty}^{+\infty}$ is an almost periodic sequence.
Since $c$ is given by (\ref{fpo}), it is the only bounded solution
of (\ref{disc}) and satisfies (\ref{estim-c}).  $\Box$


Then the following result is obtained.

\begin{teo}
\label{pre-conti} Consider $\theta$ defined in (\ref{theta}). Assume
that hypotheses (H1), (H3) and (H6) are satisfied. Then, DEPCAG
(\ref{depca1}) has a unique almost periodic solution. Moreover,
\begin{equation}
\label{estim-y} |y|_{\infty} \leq K_3|f|_{\infty},
\end{equation}
where  $\displaystyle
K_3=\left[\sqrt{q}K_0(1+|B|_{\infty}\theta)\frac{2K}{1-\rho}+1\right]\sqrt{q}K_0\theta$.
\end{teo}

\iz{Proof:} Let $\varepsilon>0$. By Lemma \ref{master}, there is
$\tau\in T(A,\varepsilon) \cap T(B,\varepsilon) \cap
T(f,\varepsilon)$ and $p \in P_{\varepsilon} \cap T(c,\varepsilon)$.
Let $y$ be the solution of (\ref{depca1}). Fix $t \in \mathbb{R}$
and let $n \in \mathbb{Z}$ such that $t \in J_n$. Then,
\[
\begin{array}{rcl}
y(t+\tau)-y(t)&=&[X(t+\tau,t_{n+p})-X(t,t_n)]c(n+p)\\
\\
&+&X(t,t_n)[c(n+p)-c(n)]\\
\\
&+&\int_{t_{n+p}-\tau}^t[X(t+\tau,u+\tau)-X(t,u)]B(u+\tau)du\cdot c(n+p)\\
\\
&+&\int_{t_{n+p}-\tau}^tX(t,u)B(u+\tau)du\cdot[c(n+p)-c(n)]\\
\\
&+&\int_{t_{n+p}-\tau}^tX(t,u)[B(u+\tau)-B(u)]du\cdot c(n)\\
\\
&+&\int_{t_{n+p}-\tau}^{t_n}X(t,u)B(u)du\cdot c(n)\\
\\
&+&\int_{t_{n+p}-\tau}^t[X(t+\tau,u+\tau)-X(t,u)]f(u+\tau)du\\
\\
&+&\int_{t_{n+p}-\tau}^tX(t,u)[f(u+\tau)-f(u)]du\\
\\
&+&\int_{t_{n+p}-\tau}^{t_n}X(t,u)f(u)du\\
\\
\end{array}
\]
So, by Lemmas \ref{conti} and \ref{zoom}, there is $K'>0$ large
enough such that
$|y(t+\tau)-y(t)| \leq \varepsilon K'$ for all $t \in \mathbb{R}$.
Since $\tau>0$ was taken arbitrarily in $T(A,\varepsilon) \cap
T(B,\varepsilon) \cap T(f,\varepsilon)$, this set is contained in
$T(x,\varepsilon K')$. By Lema \ref{master}, $T(x,\varepsilon K')$
is relatively dense. Since $\varepsilon>0$ was taken arbitrarily,
$y$ is an almost periodic solution of (\ref{depca1}). From
(\ref{pre-22}), it can be noticed that $y$ is the unique bounded
solution of DEPCAG (\ref{depca1}). So, $y$ is the unique almost
periodic solution of DEPCAG (\ref{depca1}).

Since $Z(t,s)$ is bounded and the relations (\ref{theta}),
(\ref{pre-22}) and (\ref{estim-c}) are satisfied, we have inequality (\ref{estim-y}). $\Box$
\section{The Nonlinear DEPCAG (\ref{depca2})}
To study the existence of an almost periodic solution of the DEPCAG
(\ref{depca2}), it is reminded that
$W \subseteq (\mathbb{C}^q)^{\ell}$ is not empty and the set
\[
T(F,\varepsilon,W)=\{\tau \in \mathbb{R}:|F(t+\tau,w)-F(t,w)| \leq
\varepsilon,\;\mbox{for all}\;(t,w) \in \mathbb{R}\times W\}
\]
is relatively dense for all $\varepsilon>0$.
\begin{lema} Let $y: \mathbb{R} \to \mathbb{C}^q$ an almost periodic function. Assume that (H2) is satisfied and  $F$ satisfies (H4). Then $F(t,y_\gamma(t))$ satisfies (H3).
\end{lema}
{\iz Proof:} Let $\varepsilon>0$ and $\tau \in T(y,\varepsilon) \cap
T(F,\varepsilon,W)$.
Since $y$ is almost periodic, it is uniformly continuous. So, there
is $\delta>0$ such that $s,t \in \mathbb{R}: |s-t| \leq \delta$
implies that $|y(t)-y(s)| \leq \varepsilon$. Since $P_{\tau}(\delta)
\neq \phi$,
$\left|\gamma^{p_j}(t+\tau)-\left(\gamma^{p_j}(t)+\tau\right)\right|\leq
\delta$, for $j=1,\ldots,\ell$. Moreover,
\[
\begin{array}{rcl}
|F(t+\tau,y_{\gamma}(t+\tau))-F(t,y_{\gamma}(t))| &\leq&|F(t+\tau,y_{\gamma}(t+\tau))-F(t,y_{\gamma}(t+\tau))|\\
\\
&+&|F(t,y_{\gamma}(t+\tau))-F(t,y_{\gamma}(t))|\\
\\
&\leq&\varepsilon+L\ell\varepsilon.
\end{array}
\]
Since $\varepsilon>0$ was taken arbitrarily, $F(t,y_{\gamma}(t))$ satisfies (H3).
\begin{flushright}
$\Box$
\end{flushright}
Then, the following result is obtained.
\begin{teo}
\label{apdepca2} Assume that (H1), (H2) and (H6) hold.  Assume that
$F$ satisfies (H4). If
\begin{equation}
\label{smallap1} 2\frac{KL\ell}{1-\rho}<1,
\end{equation}
then equation (\ref{depca2}) has an almost periodic solution.
\end{teo}
\iz{Proof:}  Let
\begin{equation}
\label{fpox} ({\cal T}c)(n)=\sum_{k=-\infty}^{+\infty}{\cal
G}(n,k)h(k,\hat{c}(k)),
\end{equation}
where $h(n,\hat{c}(n))=\int_{t_n}^{t_{n+1}}
X(t_{n+1},s)F(s,\hat{c}(n))ds$ and ${\cal G}(n,k)$ is given in
(\ref{green}) and $\hat{c}(n)=(c(n-p_1),\ldots,c(n-p_{\ell}))$.

If $c$ is a fixed point of the operator defined by (\ref{fpox}) then
$c$ is solution of the difference equation
\begin{equation}
\label{discx} c(n+1)=H(n)c(n)+h(n,\hat{c}(n)).
\end{equation}

If $c$ is almost periodic then $h(n,\hat{c}(n))$ is almost periodic.
In that case, ${\cal T}c$ is almost periodic. So, $\displaystyle
{\cal T}\left({\cal A}{\cal P}(\mathbb{Z},\mathbb{C}^q)\right)
\subseteq  {\cal A}{\cal P}(\mathbb{Z},\mathbb{C}^q)$. Moreover,
\[
|({\cal T}c_1)(n)-({\cal T}c_2)(n)|\leq
2\frac{KL\ell}{1-\rho}|c_1-c_2|_{\infty}.
\]

If (\ref{smallap1}) holds, ${\cal T}:{\cal A}{\cal
P}(\mathbb{Z},\mathbb{C}^q) \to {\cal A}{\cal
P}(\mathbb{Z},\mathbb{C}^q)$ is a contracting mapping. By the Banach
fixed point theorem, there is $c \in {\cal A}{\cal
P}(\mathbb{Z},\mathbb{C}^q)$ a unique fixed point for ${\cal T}$.

Therefore, equation (\ref{discx}) has an almost periodic solution
$c$. By Theorem \ref{pre-conti}, it can be constructed a solution
$y$ of (\ref{depca2}) which is almost periodic. $\Box$

From now on we will be devoted,  to the
exponential stability of the almost periodic solution of
(\ref{depca2}) whose existence was proved in the previous section. So, we say what we will understand by exponential stability.

Assume that $p_j>0$ for $j=1,\ldots,\ell$. Let $\displaystyle
p=\max_{j=1,\ldots,\ell} p_j$.

A solution $y$ of the DEPCAG (\ref{depca2}), is {\it exponentially
stable} as $t \to +\infty$ if there is $\alpha \in ]0,1[$ such that given
$\varepsilon>0$, there exists $\delta>0$ such that  $\tilde{y}=\tilde{y}(t)$ is a solution of
(\ref{depca2}) defined for $t \geq t_{0}$  then
\[
\max_{j=0,1,...,p}|y(t_{-j})-\tilde{y}(t_{-j})|
\leq \delta
\]
implies
\begin{equation}
\label{eadef} |\tilde{y}(t)-y(t)| \leq \varepsilon
\alpha^t,\;\mbox{for all}\;t \geq t_{0}.
\end{equation}
This kind of stability is in the half axis although the solution
being exponentially stable is defined on the whole axis.

The recent definition is independent on the choice of $t_0$. Any other value could be chosen.

Let $\Phi(n,k)=\Phi(n)\Phi(k)^{-1}$, for all $(n,k)\in
\mathbb{Z}^2$. Assume  that the difference system (\ref{disdich}) is
{\it exponentially stable} as $n \to +\infty$, i.e., assume that there are  positive constants
$\rho, K$ with $\rho<1$  and $K \geq 1$ such that
\begin{equation}
\label{est} |\Phi(n,k+1)| \leq K\rho^{n-k},
\end{equation}
for all $n,k \in \mathbb{Z}: n \geq k$.


By Theorem \ref{apdepca2} and the exponential stability, the
condition
\begin{equation}
\label{smallap2} \frac{KL\ell}{1-\rho}<1,
\end{equation}
insures the existence of a unique almost periodic solution $y=y(t)$
of DEPCAG (\ref{depca2}) defined for all $t \in \mathbb{R}$.

For DEPCAG (\ref{homo2}), notice that an exponential stability for
(\ref{disdich}) implies a direct notion on exponential stability on $Z(t,s)$. In fact, from
(\ref{productoria}) and (\ref{est}), it is obtained, for $n>k$,
$t\in J_n$ and $s \in ]t_k,t_{k+1}]$, that
\[
|Z(t,s)| \leq K_4 \rho^{n-k},
\]
where $K_4=K\sqrt{q}K_0^2[1+\sqrt{q}K_0|B|_{\infty}\theta]^2$ and
$\theta$ is given in (\ref{theta}). Since $t-s \leq t_{n+1}-t_k \leq
\theta (n-k+2)$,
\[
|Z(t,s)| \leq K_4\rho^{-2}\rho^{\frac{t-s}{\theta}}.
\]

If $\eta_0,\eta_1,\ldots,\eta_{p} \in
\mathbb{C}^q$, it is not hard to see that the difference system
(\ref{discx}) has a solution $\tilde{c}=\tilde{c}(n)$ defined for $n
\geq 0$ with the initial conditions $\tilde{c}(-j)=\eta_j \in
\mathbb{C}^q$ for $j=0,1,\ldots,p$.

Let
\begin{equation}
\label{22}
\begin{array}{rcl}
\tilde{y}(t)&=&Z_{k(t)}(t)\cdot
\left(\Phi(n,0)\tilde{c}(0)\right.\\
\\
&+&\left.\sum_{k=0}^{n-1}\Phi(n,k+1)\int_{t_k}^{t_{k+1}} X(t_{k+1},u)F(u,\tilde{c}(k-p_1),\ldots,\tilde{c}(k-p_{\ell}))du\right)\\
\\
&+&\int_{\gamma^0(t)}^t
X(t,u)F(u,\tilde{c}(n-p_1),\ldots,\tilde{c}(n-p_{\ell}))du,
\end{array}
\end{equation}
where   $t \geq t_0$. Then, $\tilde{y}=\tilde{y}(t)$
is the unique solution of (\ref{depca2}) with $t \geq t_{0}$ and
fixed initial conditions $\tilde{y}(t_{-j})=\eta_j$ for
$j=0,1,\ldots,p$.

Then, the following result is given.

\begin{teo} Assume that (H1), (H2), (H4)  hold and that the difference system (\ref{disdich}) has an exponential stability as $n \to +\infty$. Assume that (\ref{est}) and (\ref{smallap2}) hold.   If $y$ is the almost periodic solution  of (\ref{depca2}) and $\tilde{y}$ is solution of (\ref{depca2}) for $t \geq t_{0}$ with initial conditions $\tilde{y}(t_{-j})=\eta_j$ for $j=0,1,\ldots,p$, then there is $\tilde{K}>0$ such that
\begin{equation}
\label{es} |y(t)-\tilde{y}(t)| \leq \tilde{K}
\left({\rho}(1+KL\ell\rho^{-p})\right)^{n}\max_{j=0,1,...,p}|c(-j)-\eta_j|,
\end{equation}
where $t \geq t_0$. Hence, if
\begin{equation}
\label{smallap3} \frac{KL\ell}{1-\rho}<\rho^{p-1}
\end{equation}
then $y$ is exponentially stable.
\end{teo}

\iz{Proof:} Consider that $c(n)=y(t_n)$ and
$\tilde{c}(n)=\tilde{y}(t_n)$  for all integer $n \geq n_0$.  Let
$u(n)=c(n)-\tilde{c}(n)$ for all $n \in \mathbb{Z}$. Then, for $n_0
\in \mathbb{Z}$,
\[
\begin{array}{rcl}
|u(n)|&\leq&|\Phi(n,0)||u(0)|+\sum_{k=0}^{n-1}|\Phi(n,k+1)||F(k,\hat{{c}})(n)-F(k,\hat{\tilde{c}}(n))|\\
\\
&\leq&K\rho^{n}|u(0)|+
KL\sum_{k=0}^{n-1}\rho^{n-k}\sum_{j=1}^{\ell}|u(k-p_j)|.
\end{array}
\]

Let $\displaystyle \omega(n)=\rho^{-n}\sum_{j=1}^{\ell}|u(n-p_j)|$
and  $v(n)=\rho^{-n}|u(n)|$. Then
\[
\rho^{-n}u(n)\leq K|u(0)|+
KL\sum_{k=0}^{n-1}\omega(k)
\]

Notice that $\displaystyle \omega(n)=\sum_{j=1}^{\ell}
\rho^{-p_j}\rho^{(-n-p_j)}|u(n-p_j)|\leq \rho^{-p}\sum_{j=1}^{\ell}
v(n-p_j)$. For $n \geq 0$,
\[
v(n) \leq Kv(0)+KL\sum_{k=0}^{n-1} \omega(k).
\]

Let $\displaystyle z_n=\max\{|v(m)|:m=-p,-p+1,\ldots,n\}$.
Then, $\omega(n) \leq \rho^{-p}\ell z_n$, for all $n \geq 0$.
Hence,
\[
v(n) \leq Kv(0)+KL\rho^{-p}\ell\sum_{k=0}^{n-1} z_k.
\]

Let $m_n \in \{-p,n-p+1,\ldots,n\}$ such that $z_n=v(m_n)$.


If $m_n \geq 0$, then $\displaystyle z_n \leq
Kv(0)+KL\ell\rho^{-p}\sum_{k=0}^{m_n-1}z_k$.
Hence,
\[
z_n\leq Kv(0)+KL\ell\rho^{-p}\sum_{k=0}^{n-1} z_k.
\]
If $m_n < 0$, then there is $j_0 \in \{1,\ldots,p\}$ such that
$m_n=n-j_0$. Since $K \geq 1$, $z_n \leq K z_{0}$. So,
\[
z_n\leq Kz_{0}+KL\ell\rho^{-p}\sum_{k=0}^{n-1} z_k,
\]
for all $n \geq 0$.
By Gronwall's inequality,
\[
z_n \leq (1+KL\ell\rho^{-p})^{n}z_{0}.
\]
So, for all $n \geq 0$,
\begin{equation}
\label{as1} \displaystyle |c(n)-\tilde{c}(n)| \leq
K\rho^{n}(1+KL\ell\rho^{-p})^{n}\max_{j=0,1,\ldots,p}|c(-j)-\tilde{c}(-j)|.
\end{equation}

By Lemma \ref{conti}, there is a positive constant $K_0$ such that
$|X(t,u)| \leq \sqrt{q}K_0$, for all $u \in J_{k(t)}$ and $|Z(t,\gamma^0(t))| \leq
\sqrt{q}K_0(1+\sqrt{q}K_0|B|_{\infty}\theta)$ for all $t \geq t_0$.
By relation (\ref{22}), for $t \geq t_0$, there is
a positive constant $K'$ such that
\begin{equation}
\label{as2} |y(t)-\tilde{y}(t)| \leq K'|c(n)-\tilde{c}(n)|.
\end{equation}

The recent inequality show a Lipschitz continuous relation $\tilde{c} \mapsto \tilde{y}$.

By combining (\ref{as1}) and (\ref{as2}), this result is proved with
$\tilde{K}=K'K$.

Notice that (\ref{smallap3}) implies (\ref{smallap2}). Then Theorem
\ref{apdepca2} insures the existence of the unique almost periodic
solution of (\ref{depca2}) which is exponentially stable.  In fact,
let $\alpha={\rho}(1+KL\ell\rho^{-p})$. By (\ref{smallap3}),
$\alpha<1$. For $\varepsilon>0$ consider
$\delta=\frac{\varepsilon}{\tilde{K}}$. By (\ref{es}),
(\ref{eadef}) is satisfied and $y$ is
exponentially stable.

\begin{flushright}
$\Box$
\end{flushright}
In the last theorem, the condition (\ref{smallap3}) is simple and
slightly stronger than the condition of existence (\ref{smallap2}).

\section{Examples}
\subsection{Exponential Dichotomy}
It is not obvious to extend the exponential dichotomy for the
difference equation (\ref{disdich}) for the DEPCAG (\ref{homo2}). We
could consider an intuitively direct  definition given by the
existence of a projection $\Pi_*$ and positive constants $M$ and
$\alpha$ such that
\begin{equation}
\label{intento}
\begin{array}{rcl}
|Z(t,t_0)\Pi_*Z(s,t_0)^{-1}| \leq Me^{-\alpha(t-s)},&\mbox{if}& t \geq s\\
\\
|Z(t,t_0)(I-\Pi_*)Z(s,t_0)^{-1}| \leq Me^{\alpha(t-s)},&\mbox{if}&
t\leq s.
\end{array}
\end{equation}
However, if we take $A(t)=0$,
$B(t)=\mbox{diag}(\lambda_0(t),\lambda_1(t))$,
$\displaystyle\lambda_0(t)=-\frac{2}{\pi}+\sin(2\pi t)$,
$\lambda_1(t)=-\lambda_0(t)$,  $t_n=n$ for all $n \in \mathbb{Z}$,
$\displaystyle
\int_{n}^{n+\delta}\lambda_0(\xi)d\xi=-\frac{1}{2\pi}\left(4\delta-1+\cos(2\pi\delta)\right)$
and $\displaystyle
\int_{n}^{n+\delta}\lambda_1(\xi)d\xi=\frac{1}{2\pi}\left(4\delta-1+\cos(2\pi\delta)\right)$
for all $\displaystyle \delta \in [0,1]$ then the exponential
dichotomy on the difference equation (\ref{disdich}) which can be
written  as (\ref{edic}) is satisfied for $\Pi=\mbox{diag}(1,0)$ but
there is no $\Pi_*$ such that condition (\ref{intento}) is
satisfied.

Notice that a dichotomy condition on the ordinary differential
equation (\ref{homo}) implies an exponential dichotomy on the
difference equation (\ref{disdich}) \cite[Proposition 2]{Papas94}
when $|B(t)|$ is small enough and $y_{\gamma}(t)=y([t])$. However,
an exponential dichotomy for the difference equation on
(\ref{disdich}) is not a necessary condition for an exponential
dichotomy for the ordinary differential system (\ref{homo}). In
fact, let's consider $t_n=n$, $A(t)=0$ and
$B(t)=\mbox{diag}\left(-\frac{3}{2},\frac{1}{2}\right)$. Then the
exponential dichotomy for difference system (\ref{disdich}) is
satisfied, with no exponential dichotomy for the ordinary
differential system (\ref{homo}).

\subsection{Constant Coefficients}

Assume that in DEPCAG (\ref{depca2}), $A(t)=A_0$ and $B(t)=B_0$ are
constants matrices and $F(t,\cdot)$ is almost periodic. Then DEPCAG
(\ref{depca2}) becomes
\begin{equation}
\label{depca3}
y'(t)=A_0y(t)+B_0y(\gamma^{0}(t))+F(t,y_{\gamma}(t)),\;t\in
\mathbb{R}.
\end{equation}
Assume that $t_{n+1}-t_n=\nu$ constant, $A_0$ and
\[
H(n)=H_0=e^{\nu
A_0}\left[I+A_0^{-1}(I-e^{-\nu A_0})B_0\right]
\]
are invertible.

By using $\sigma(H_0)$ as the usual notation for the spectrum of the
matrix $H_0$, it can be said that:

If $\sigma(H_0) \cap \{ z \in \mathbb{C}:|z|=1\}$ is the empty set
and $L$ in (\ref{lp}) satisfies (\ref{smallap2}), then DEPCAG
(\ref{depca3}) has an almost periodic solution. In particular, it is
obtained when the elements of $\sigma(A_0)$ have non zero real part
and $|B_0|$ is small enough.

If $ \sigma(H_0) \subseteq \{ z \in \mathbb{C}:|z|<1\}$ and $L$ in
(\ref{lp}) satisfies (\ref{smallap3}), then DEPCAG (\ref{depca3})
has an almost periodic solution which is exponentially stable.  In
particular, it is obtained when the elements of $\sigma(A_0)$ have
negative real part and $|B_0|$ is small enough.

If $A_0=0$, $H(n)=I+\nu B_0$ invertible. If $\sigma(B_0) \subseteq
\{ z \in \mathbb{C}:|z|<1/\nu\}$ is the empty set and $L$ in
(\ref{lp}) satisfies (\ref{smallap2}), then DEPCAG (\ref{depca3})
has an almost periodic solution. If $ \sigma(B_0) \subseteq \{ z \in
\mathbb{C}:|z|<1/\nu\}$ and $K$ in (\ref{lp}) satisfies
(\ref{smallap3}), then DEPCAG (\ref{depca3}) has an almost periodic
solution which is exponentially stable. We can notice that it
behaves as a difference equation.


\begin{thebibliography}{10}
\bibitem{AL} E. Ait Dads and L. Lhachimi. Pseudo almost periodic solutions for equation with piecewise constant argument. {\it J. of Math. Analysis and Applications.} {\bf 371} (2010), 842--854.
\bibitem{A} M. U. Akhmet. Integral manifolds of differential equations with piecewise
constant argument of generalized type. {\it Nonlinear Analysis TMA}
{\bf 66} (2) (2007), 367--383.
\bibitem{A02} M.U. Akhmet. Almost periodic solutions of differential equations with piecewise constant argument of generalized type. {\it Nonlinear Analysis: Hybrid Systems.} {\bf 2} (2008) 456--467.
\bibitem{A03} M.U. Akhmet. Stability of differential equations with piecewise constant arguments of generalized type. {\it Nonlinear Analysis: TMA} {\bf 68} (2008) 794--803.
\bibitem{A04}  M.U. Akhmet, C. Buyukadali. Differential equations with state-dependent piecewise constant argument, {\it Nonlinear Analysis: TMA} {\bf 72} (2010) 4200--4210.
\bibitem{Be} A. S. Besicovitch. {\it Almost Periodic Functions.} Dover Publications, Inc., New York, 1955.
\bibitem{Bochner} S. Bochner. Beitrage zur Theorie der
fastperiodische Funktionen. I: Funktionen einer Variablen. {\it
Math. Ann.} {\bf 96} (1927) 119--147.
\bibitem{Bo01} H. Bohr. Zur Theorie der fastperiodischen Funktionen I. {\it Acta Math.} {\bf 45} (1925)  29--127.
\bibitem{Bo02} H. Bohr. {\it Almost Periodic Functions.} Chelsea Publishing Company, New York, N.Y., 1947.
\bibitem{X2005} J. Cao, M. Lin and Y. Xia. The existence of almost periodic solutions of certain perturbation systems. {\it J. of Math. Analysis and Applications.} {\bf 310} (2005), 81--96.
\bibitem{X2005-1} J. Cao and Y. Xia. Almost periodicity in an ecological model with $M$-predators and $N$-preys by ``pure-delay type'' system. {\it Nonlinear Dynamics} {\bf 39} (2005), 275--304.
\bibitem{Ch2011} K. Chiu. Stability of oscillatory solutions of differential equations with a general piecewise constant argument. {\it Electron. J. Qual. Theory Differ. Equ.} (2011) {\bf 88}, 15 pp.
\bibitem{ChP2010} K. Chiu and M. Pinto. Periodic solutions of differential equations with a general piecewise constant argument and applications. {\it Electron. J. Qual. Theory Differ. Equ.} (2010) {\bf 46}, 19 pp.
\bibitem{CoWi84} K. L. Cooke and J. Wiener. Retarded differential equations with
piecewise constant delays. {\it J. of Math. Analysis and
Applications.} {\bf 99} (1984), 265--297.
\bibitem{Cop78} W. A. Coppel. {\it Dichotomies in Stability Theory.} Lecture Notes in Mathematics, Vol. 629. Springer-Verlag, Berlin-New York, 1978.
\bibitem{CP} C. Cuevas and M. Pinto. Existence and uniqueness of pseudo almost periodic solutions of semilinear Cauchy problems with non dense domain. {\it Nonlinear Analysis TMA.} {\bf 45} (8) (2001)  73--83.
\bibitem{F5} A.M. Fink. {\it Almost Periodic Differential Equations.}  Lecture Notes in Mathematics {\bf 377} Springer, Berlin, 1974.
\bibitem{X2007}  M. Han, Z. Huang and Y. Xia. Existence of almost periodic solutions for forced perturbed systems with piecewise constant argument. {\it J. of Math. Analysis and Applications.} {\bf 333} (2007), 798--816.
\bibitem{YH16} J. Hong and R. Yuan. The existence of almost periodic solutions for a class of differential equations with
piecewise constant argument. {\it Nonlinear Analysis TMA.} {\bf 28}
(8) (1997)\\ 1439--1450.
\bibitem{Papas94} G. Papaschinopoulos. Exponential dichotomy, topological equivalence and structural stability for differential equations with piecewise constant argument. {\it Analysis} {\bf 14} n 2-3  (1994) 239--247.
\bibitem{SP} N.A. Perestyuk and A.M. Samoilenko. {\it Impulsive Differential Equations.} World Scientific, 1995.
\bibitem{Pi2009} M. Pinto. Asymptotic equivalence of nonlinear and quasi linear differential equations with piecewise constant arguments. {\it Mathematical and Computer Modelling.} {\bf 49} (2009) 1750--1758.
\bibitem{Pi2010} M. Pinto. Dichotomy and existence of periodic solutions of quasilinear functional differential equations. {\it Nonlinear Analysis TMA.} {\bf 72} (3-4) (2010) 1227--1234.
\bibitem{Pi2010JDEA}Pinto, Manuel. Cauchy and Green matrices type and stability in alternately advanced and delayed differential systems.{\it J. of Difference Equations and Applications}. {\bf 17} (2011), no. 2, 235--254.

\bibitem{PAP} M. Pinto. Pseudo-almost periodic solutions of neutral integral and differential equations with applications. {\it Nonlinear Analysis TMA.} {\bf 72} (2010) 4377--4383.
\bibitem{PR} M. Pinto and G. Robledo. Existence and stability of almost periodic solutions in impulsive neural network models. {\it Applied Mathematics and Computation} {\bf 217} (2010), no. 8, 4167--4177.
\bibitem{ShWi83} S. M. Shah and J. Wiener. Advanced differential equations with piecewise
constant argument deviations. {\it Internat. J. Math. and Math.
Sci.} {\bf 6} n 4 (1983) 671--703.
\bibitem{R3} N. Van Minh.  Almost periodic solutions of $C$-well-posed evolution equations. {\it Math. J. Okayama Univ.} {\bf 48} (2006), 145--157.
\bibitem{Wi83} J. Wiener. Differential equations with piecewise constant delays. {\it Trends in theory and practice of nonlinear differential equations (Arlington, Tex., 1982)},  {\it Lecture Notes in Pure and Appl. Math.} {\bf 90}  Dekker, New York, (1984) 547--552.
\bibitem{Wi93} J. Wiener. {\it Generalized Solutions of Functional Differential
Equations.} World Scientific, 1993.
\bibitem{Y} R Yuan. The existence of almost periodic solutions of retarded differential equations with piecewise constant argument. {\it Nonlinear Analysis TMA.} {\bf 48} (7) (2002), 1013--1032.
\bibitem{Z} C. Zhang. {\it Almost Periodic Type Functions and Ergodicity.} Science Press, Beijing; Kluwer Academic Publishers, Dordrecht, 2003.
\end{thebibliography}
\end{document}